\documentclass[11pt]{amsart}
\usepackage{amssymb,amsmath,hyperref,breqn}

\newtheorem{theorem}{Theorem}

\newtheorem{conjecture}[theorem]{Conjecture}
\newtheorem{corollary}[theorem]{Corollary}
\newtheorem{definition}[theorem]{Definition}
\newtheorem{example}[theorem]{Example}

\newtheorem{fact}[theorem]{Fact}
\newtheorem{lemma}[theorem]{Lemma}
\newtheorem{problem}[theorem]{Problem}
\newtheorem{proposition}[theorem]{Proposition}
\newtheorem{question}[theorem]{Question}
\newtheorem{remark}[theorem]{Remark}

\newcommand{\bcon}{\begin{conjecture}}
\newcommand{\econ}{\end{conjecture}}
\newcommand{\bcor}{\begin{corollary}}
\newcommand{\ecor}{\end{corollary}}
\newcommand{\bdf}{\begin{definition}}
\newcommand{\edf}{\end{definition}}
\newcommand{\beq}{\begin{equation}}
\newcommand{\eeq}{\end{equation}}
\newcommand{\bexa}{\begin{example}}
\newcommand{\eexa}{\end{example}}
\newcommand{\bfac}{\begin{fact}}
\newcommand{\efac}{\end{fact}}
\newcommand{\blem}{\begin{lemma}}
\newcommand{\elem}{\end{lemma}}
\newcommand{\bprb}{\begin{problem}}
\newcommand{\eprb}{\end{problem}}
\newcommand{\bpro}{\begin{proposition}}
\newcommand{\epro}{\end{proposition}}

\newcommand{\bque}{\begin{question}}
\newcommand{\eque}{\end{question}}
\newcommand{\brem}{\begin{remark}}
\newcommand{\erem}{\end{remark}}
\newcommand{\bthm}{\begin{theorem}}
\newcommand{\ethm}{\end{theorem}}
\newcommand{\bmat}{\begin{pmatrix}}
\newcommand{\emat}{\end{pmatrix}}

\newcommand{\be}{\begin}
\newcommand{\en}{\end}
\newcommand{\bpr}{\begin{proof}}
\newcommand{\epr}{\end{proof}}

\newcommand{\lb}{\label}

\newcommand{\comment}[1]{\,}

\newcommand{\cal}{\mathcal}

\newcommand{\N}{\mathbb N}
\newcommand{\Z}{\mathbb Z}

\newcommand{\C}{\mathbb C}

\newcommand{\ve}{\varepsilon}
\newcommand{\subsetn}{\varsubsetneq}

\title{$SO(2n,\C)$-character varieties are not varieties of characters}

\author{Adam S. Sikora}

\thanks{The author acknowledges support from U.S. National Science Foundation grants DMS 1107452, 1107263, 1107367 "RNMS: GEometric structures And Representation varieties" (the GEAR Network).}


\begin{document}

\thispagestyle{empty}

\begin{abstract}
We prove that the coordinate rings of $SO(2n,\C)$-character varieties are not generated by trace functions nor generalized trace functions for $n\geq 2$ and all groups $\Gamma$ of corank $\geq 2.$ 
Furthermore, we give examples of non-conjugate completely reducible representations undistinguishable by generalized trace functions. Hence,  $SO(2n,\C)$-character varieties are not varieties of characters!

However, we also prove that any generic $SO(2n,\C)$-representation of a free group can be distinguished from all non-equivalent representations by trace functions and by a single generalized trace function.
\end{abstract}

\address{244 Math Bldg, University at Buffalo, SUNY, Buffalo, NY 14260}
\email{asikora@buffalo.edu}

\maketitle

\pagestyle{myheadings}

%
%

For a complex, reductive algebraic group $G$, and a finitely generated discrete group $\Gamma$, the $G$-character variety of $\Gamma$, $X_G(\Gamma),$ 
is the categorical quotient of the representation variety, $Rep(\Gamma,G)$, by the action of $G$ by conjugation, cf. \cite{LM, S-char} and the references within\footnote{Throughout the paper, the field of complex numbers can be replaced an arbitrary algebraically closed field of characteristic zero.}. (In this paper $Rep(\Gamma,G)$ and $X_G(\Gamma)$ are affine algebraic sets rather than possibly non-reduced schemes of \cite{LM, S-char}).

If $G$ is a matrix group, then every $\gamma\in \Gamma$ defines a trace function 
$$\tau_\gamma: X_G(\Gamma) \to \C, \quad \tau_\gamma([\rho])=tr\rho(\gamma).$$ 
The $G$-trace algebra of $\Gamma$, denoted by ${\cal T}_{G}(\Gamma)$, is a subalgebra  of $\C[X_{G}(\Gamma)]$ generated by $\tau_\gamma$, for all $\gamma\in \Gamma.$

For $G=SL(n,\C), Sp(n,\C)$ (symplectic groups), and $SO(2n+1,\C)$ (odd special orthogonal groups) the coordinate rings, $\C[X_G(\Gamma)],$ are generated by trace functions, \cite{S-gen, FL}, and although one might expect that for other groups $G$ as well, it is not the case for $G=SO(2n,\C)$. Indeed, if $\rho:\Gamma\to SO(2n,\C)$ is irreducible and $\rho'$ is obtained from $\rho$ by a conjugation by a matrix $M\in O(2n,\C),$ $det M=-1,$ then $\rho$ and $\rho'$ are indistinguishable by trace functions even though it is not difficult to see that they are not $SO(2n,\C)$-conjugate and, furthermore, they are distinct in $X_{SO(2n,\C)}(\Gamma),$
cf. \cite[Sec. 4]{S-quo}.

Pursuing better conjugacy invariants of representations, one is naturally lead to the notion of generalized trace functions, which (as an added bonus) are defined irrespectively of specific realizations of $G$ as matrix groups. More specifically, for any $\gamma\in \Gamma$ and any (finite dimensional) representation $\phi$ of $G$ there is the generalized trace function 
$$\tau_{\phi, \gamma}:X_G(\Gamma) \to \C,\quad  \tau_{\phi,\gamma}([\rho])=tr\phi\rho(\gamma).$$  
The full $G$-trace algebra of $\Gamma$, ${\cal FT}_{G}(\Gamma),$ is a subalgebra of $\C[X_{G}(\Gamma)]$
generated by $\tau_{\phi,\gamma}$, for all $\gamma\in \Gamma$ and all $\phi$, cf. \cite{S-gen}.

We proved in \cite{S-quo} however that even these functions do not generate $\C[X_{G}(\Gamma)]$ in general. Specifically, for any group $\Gamma$ of corank $\geq 2$ (ie. having the free group on two generators, $F_2$, as a quotient) ${\cal FT}_{SO(4,\C)}(\Gamma)$ is a proper subalgebra of $\C[X_{SO(4,\C)}(\Gamma)]$.
The main result of this paper is a generalization of this statement to even orthogonal groups of higher rank.

\bthm \lb{main}
\beq\lb{prop_subalg}
{\cal FT}_{SO(2n,\C)}(\Gamma)\subsetn \C[X_{SO(2n,\C)}(\Gamma)]
\eeq
for every group $\Gamma$ of corank $\geq 2$ and all $n\geq 2.$
\ethm

Let $\imath: SO(4,\C)\to SO(2n,\C)$ for $n>2$ be the ``obvious'' embedding (extending $4\times 4$ matrices by the $(2n-4)\times (2n-4)$ identity matrix). The above result would be an easy consequence of our earlier result for $SO(4,\C)$ if the induced homomorphism  
$$\imath_*: \C[X_{SO(2n,\C)}(\Gamma)] \to \C[X_{SO(4,\C)}(\Gamma)]$$ 
was onto. However, one can prove that the image of $\imath_*$ is ${\cal T}_{SO(4,\C)}(\Gamma),$ cf. Corollary \ref{obvious_emb}.

It remains an open question whether the proper inclusion (\ref{prop_subalg}) holds if $F_2$ is replaced by some of its quotients, like $\Z^2.$ 
Our proof in \cite{S-quo} of the result for the free group $\Gamma=F_2$ and $n=2$ relies on the fact that $\C[X_{SL(2,\C)}(\Gamma)]$ has an $\N$-grading, which is not the case for $\Gamma=\Z^2$ and, therefore, it does not generalize to that case.  

It is worth mentioning though that we can prove that 
$${\cal FT}_{SO(4,\C)}(\Gamma)=\C[X_{SO(4,\C)}(\Gamma)]$$
for $\Gamma=\Z_p^2,$ for $p=2,3,4,5$ by a direct computation.
(We use the topologist's notation, $\Z_p,$ for the cyclic group of order $p$.)

Note that the proper inclusion (\ref{prop_subalg}) does not imply that the generalized trace functions fail to distinguish non-equivalent $SO(2n,\C)$-representations of $\Gamma.$ 
To see that consider as an example representations into $GL(n,\C)$. Their equivalence classes in $X_{GL(n,\C)}(\Gamma)$ are distinguished by trace functions, even though
$${\cal T}_{GL(n,\C)}(\Gamma)\subsetn \C[X_{GL(n,\C)}(\Gamma)],$$
since it is easy to show that
$$f_\gamma([\rho])=det(\rho(\gamma))^{-1}$$ 
is a regular function on $X_{GL(n,\C)}(\Gamma)$ which is not a polynomial in trace functions
for example for $\Gamma=\Z$ and $\gamma\ne 0.$


Nonetheless, we can strengthen our Theorem \ref{prop_subalg} for higher $n$ as follows.

\bthm \lb{undisting} 
For any $n\geq 7,$ $n\ne 8,$ any group $\Gamma$ projecting onto $\Z_p*\Z_q$ for $p,q>max(2n-14,16)$
has two non-conjugate irreducible representations into $SO(2n,\C)$ which are undistinguishable by the generalized trace functions.
 
These representations are related one to another through the involution $\sigma$ of \cite[Sec. 4]{S-quo}.
 \ethm
  
Since the elements of $\C[X_{SO(2n,\C)}(\Gamma)]$ distinguish points of $X_{SO(2n,\C)}(\Gamma),$ the proper inclusion (\ref{prop_subalg}) holds for all 
$\Gamma$ satisfying the assumptions of Theorem \ref{undisting}.
 
The embedding 
$${\cal FT}_{SO(2n,\C)}(\Gamma)\subset \C[X_{SO(2n,\C)}(\Gamma)]$$
and the dual map to it which we denote by $\psi,$
$$\psi: X_{SO(2n,\C)}(\Gamma) \to Spec\, {\cal FT}_{SO(2n,\C)}(\Gamma),$$
is analyzed in the next statement. 

\bpro\lb{normalization} Let $n\geq 1$ and let $\Gamma$ be a group such that
$X_{SO(2n,\C)}(\Gamma)$ is irreducible and 
$\cal T_{SO(2n,\C)}(\Gamma)\ne \cal FT_{SO(2n,\C)}(\Gamma)$.
Then\\
(1) $\psi$ 
is finite and birational. 

In particular, $\psi$ is generically $1$-$1$. It is also a normalization map if $X_{SO(2n,\C)}(\Gamma)$ is normal.\\
(2) Let $\eta$ be one of the two irreducible representations $D^\pm_n$ of \cite[\S 23.2]{FH}, with the highest weight being twice that of $\pm$-half spin representation. 
Then there exists $\gamma\in \Gamma$ such that for a generic 
$[\rho]\in X_{SO(2n,\C)}(\Gamma)$ the trace functions together with $\tau_{\eta,\gamma}$
distinguish $[\rho]$ from any other point on the $SO(2n,\C)$-character variety of $\Gamma.$
(Obviously, only a finite set of trace functions generating the trace algebra is sufficient here.)
\epro

Note that the above assumptions are satisfied by free groups of rank $\geq 2$ for $n\geq 2.$

%
\section{Proof of Theorem \ref{main}} 
%

In \cite{S-gen}, we have introduced a function $Q$ of $n$ arguments $A_1,...,A_n\in M(2n,\C)$ defined as follows\footnote{For simplicity, we skip the index $2n$ used in \cite{S-gen}.}:\\
$Q(A_1,...,A_n)=\sum_{\sigma\in S_{2n}} sn(\sigma)
 (A_{1,\sigma(1),\sigma(2)}-A_{1,\sigma(2),\sigma(1)})...\hspace*{.2in}$\vspace*{-.2in}\\
\be{equation}\label{def-Q}
\hspace*{2.6in} (A_{n,\sigma(2n-1),\sigma(2n)}-A_{n,\sigma(2n),\sigma(2n-1)}),
\en{equation}
where $A_{i,j,k}$ is the $(j,k)$-th entry of $A_i$ and $sn(\sigma)=\pm 1$ is the sign of $\sigma.$

By abuse of notation, for any $\gamma_1,...,\gamma_n\in \Gamma$ we have also defined a function $Q(\gamma_1,...,\gamma_n)$ on $Hom(\Gamma, SO(2n,\C))$ sending 
$[\rho]$ to $Q(\rho(\gamma_1),...,\rho(\gamma_n)).$ Since $Q(A_1,...,A_n)$ is invariant under simultaneous conjugation by matrices in $SO(2n,\C),$ the above factors to a function  $Q(\gamma_1,...,\gamma_n)$ on $X_{SO(2n,\C)}(\Gamma).$

Let $$SO(2n,\C)=\{A: AA^T=I,\ det(A)=1\}.$$
By \cite{S-gen}, $\C[X_{SO(2n,\C)}(\Gamma)]$ is generated by the trace functions, $\tau_\gamma,$ and the functions $Q(\gamma_1,...,\gamma_n)$, for $\gamma, \gamma_1,...,\gamma_n\in \Gamma.$

Any homomorphism 
$\pi: \Gamma\to F_2$ induces
$$\pi_*: \C[X_{SO(2n,\C)}(\Gamma)]\to \C[X_{SO(2n,\C)}(F_2)]$$
which restricts to 
$$\pi_*: {\cal FT}_{SO(2n,\C)}(\Gamma)\to {\cal FT}_{SO(2n,\C)}(F_2).$$ Furthermore, we have
$$\pi_*\tau_{\phi,\gamma}=\tau_{\phi,\pi(\gamma)}$$
and
$$\pi_* Q_n(\gamma_1,...,\gamma_n)=Q_n(\pi(\gamma_1),...,\pi(\gamma_n)).$$
Consequently,  if $\pi$ is an epimorphism then all generators of $\C[X_{SO(2n,\C)}(F_2)]$ are in the image of $\pi_*$. Hence,
it is sufficient to prove Theorem \ref{main} for $\Gamma=F_2$ only.

Therefore, assume $\Gamma=F_2$ from now on, with generators $\gamma_1, \gamma_2$. 

Note that there is a group isomorphism $\C^*\to SO(2,\C)$ sending $c\in \C^*$ to
$$D_c=\left(\begin{array}{cc} (c+c^{-1})/2 & i(c-c^{-1})/2\\
-i(c-c^{-1})/2 & (c+c^{-1})/2\\ \end{array}\right),$$
where $i^2=-1.$
For $c\in \C^*,$ consider a map $\imath_c: SO(4,\C)\to SO(2n,\C)$ sending $A\in SO(4,\C)$
to a block matrix composed of $A$ (in the top left corner) and of $n-1$ blocks
$D_c$ along the diagonal.
Using the identification of $Hom(F_2, SO(2n,\C))$ with $SO(2n,\C)^2$ through a map sending a representation $\rho$ to
$(\rho(\gamma_1),\rho(\gamma_2)),$ we have a map
$$\alpha_{c_1,c_2}\hspace*{-0.05in}: Hom(F_2, SO(4,\C))=SO(4,\C)^2\to SO(2n,\C)^2=Hom(F_2, SO(2n,\C)),$$
sending $(A_1,A_2)\in SO(4,\C)^2$ to $(\imath_{c_1}(A_1),\imath_{c_2}(A_2)).$ Note that the above map is not induced by a homomorphism from $SO(4,\C)$ to $SO(2n,\C)$ unless $c_1=c_2=1$. 

Alternatively, $\alpha_{c_1,c_2}$ can be defined as follows: Let 
$$w=(w_1,w_2): F_2\to \Z^2$$ 
be the abelianization map. Then $\alpha_{c_1,c_2}$ maps $\rho$ to $\alpha_{c_1,c_2}(\rho)$ such that 
\beq \lb{e_alphac1c2}
\alpha_{c_1,c_2}(\rho)(\gamma)=\imath_c\rho(\gamma),
\eeq 
where $c=c_1^{w_1(\gamma)}c_2^{w_2(\gamma)}$.

It is easy to see that  $\alpha_{c_1,c_2}$ factors through
$$\alpha_{c_1,c_2}: X_{SO(4,\C)}(F_2)\to X_{SO(2n,\C)}(F_2).$$

The ``obvious'' embedding $SO(4,\C)\hookrightarrow SO(2n,\C)$ induces the map
$\alpha_{1,1}.$ As mentioned in the introduction, it is natural to attempt to prove Theorem \ref{prop_subalg} by showing that the dual map 
$$\alpha_{1,1*}: \C[X_{SO(2n,\C)}(F_2)]\to \C[X_{SO(4,\C)}(F_2)]$$
is onto. As we will see in Corollary \ref{obvious_emb} below, that is unfortunately not the case. We will prove Theorem \ref{prop_subalg} by showing that $\alpha_{c_1,c_2*}$ is onto for other values of $c_1,c_2$ instead.

\bpro \lb{alpha-epi}
$\alpha_{c_1,c_2*}:\C[X_{SO(2n,\C)}(F_2)]\to \C[X_{SO(4,\C)}(F_2)]$ is onto for $c_1,c_2\ne \{-1,0,1\},$ $c_1\ne  \pm c_2^{\pm  1}.$
\epro 

The proof of this statement follows. Afterwards, we are going to prove Theorem \ref{alpha-FT-onto} asserting that 
$$\alpha_{c_1,c_2*}({\cal FT}_{SO(2n,\C)}(F_2))\subset {\cal FT}_{SO(4,\C)}(F_2).$$ These two statements together with the proper inclusion 
$${\cal FT}_{SO(4,\C)}(F_2)\subsetn \C[X_{SO(4,\C)}(F_2)]$$
proved in \cite{S-quo} imply that ${\cal FT}_{SO(2n,\C)}(F_2)$ is a proper subalgebra of\\ 
$\C[X_{SO(2n,\C)}(F_2)]$, completing the proof of Theorem \ref{main}.

We precede the proof of Proposition \ref{alpha-epi} with a few lemmas.

\blem \lb{prop_Q}
Let $A_1,...,A_n\in M(2n,\C)$. If each $A_i$ is a block matrix built of an upper left diagonal block $B_i$ of dimension $(2n-2)\times (2n-2)$ and of lower right $2\times 2$ block $C_i,$ then $Q(A_1,...,A_n)$ equals
$$\sum_{i=1}^n Q(B_1,...,\hat B_i,...,B_n)\cdot Q(C_i),$$
where $\hat\cdot$ denotes an omitted symbol.
\elem

Note that $Q\left(\be{pmatrix} a_{11} & a_{12}\\ 
 a_{21} & a_{22}\\ 
\en{pmatrix}\right)=2(a_{21}-a_{12})$ and $Q(D_c)=i(c-c^{-1}).$

\bpr By our assumptions, a summand in (\ref{def-Q}) for $\sigma\in S_n$ is not zero only if for every $i=1,...,n$ both of $\sigma(2i-1)$ and $\sigma(2i)$ belong to
$\{1,...,2n-2\}$ or to $\{2n-1,2n\}.$ Therefore, if we denote the set of $\sigma\in S_n$ such that $\sigma(2i-1),\sigma(2i)\in \{2n-1,2n\}$ by $S_n^i,$ then
$$Q(A_1,...,A_n)=\sum_{i=1}^n S_i(A_1,...,A_n),$$
where $S_i$ is the sum of terms of (\ref{def-Q}) for $\sigma\in S_n^i$. 
Now it is easy to see that formula (\ref{def-Q}) implies that 
$$S_i(A_1,...,A_n)=Q(B_1,...,\hat B_i,...,B_n)\cdot Q(C_i).$$
\epr

\bcor\lb{obvious_emb} 
The ``obvious'' embedding, $\imath_1: SO(4,\C)\hookrightarrow SO(2n,\C)$ (extending matrices by the $(2n-4)\times(2n-4)$ identity matrix), induces the map
$\alpha_{1,1*}: \C[X_{SO(2n,\C)}(\Gamma)]\to  \C[X_{SO(4,\C)}(\Gamma)]$ whose image is $\cal T_{SO(4,\C)}(\Gamma).$
\ecor

\bpr The $SO(4,\C)$ matrices embedded into $SO(2n,\C)$ as above contain a lower right block $D_1.$ Since $Q(D_1)=0,$ the above lemma implies that 
$\alpha_{1,1*}Q(\gamma_1,...,\gamma_n)=0$ 
(as a function on $X_{SO(4,\C)}(\Gamma)$) for all $\gamma_1,...,\gamma_n\in \Gamma.$

On the other hand,  since 
$$\alpha_{1,1*}\tau_\gamma=\tau_\gamma+2n-4,$$ 
the homomorphism $\alpha_{1,1*}$ maps the $SO(2n,\C)$-trace algebra onto $SO(4,\C)$-trace algebra, implying the statement. 
\epr

Denote by $Q_{k,l}(A_1,A_2)$ the value of $Q$ evaluated at $k$ matrices $A_1$ and $l$ matrices $A_2.$ For convenience, assume that $Q_{k,l}(A_1,A_2)=0$ for  $k<0$ or $l<0.$

\bcor\lb{Qkl-cor}
If $A_1,A_2$ are block matrices as in Lemma \ref{prop_Q}, then
$$Q_{k,l}(A_1,A_2)=k\cdot Q_{k-1,l}(B_1,B_2)Q(C_1)+l\cdot Q_{k,l-1}(B_1,B_2)Q(C_2).$$
\ecor

Let $Q_n(A)$ denote the value of $Q$ for all its $n$ arguments equal $A\in M(2n,\C).$

\blem \lb{Qimath}
Let $n\geq 2.$
\begin{enumerate}
\item $Q_n(\imath_{c}(A))=\frac{1}{2}[i(c-c^{-1})]^{n-2}n! Q_2(A)$, for any $A\in M(4,\C).$ \vspace*{.05in}
\item $Q_{n-1,1}(\imath_{c_1}(A_1),\imath_{c_2}(A_2))=[i(c_1-c_1^{-1})]^{n-2} (n-1)! Q_{1,1}(A_1,A_2)+$\vspace*{.05in}\\
\hspace*{1.7in}$\frac{1}{2}i(c_2-c_2^{-1}) Q_2(A_1)\sum_{k=2}^{n-1} [i(c_1-c_1^{-1})]^{k-2}k!$ for any $A_1,A_2\in M(4,\C).$
\end{enumerate}
\elem

\bpr  We have $Q(D_c)=i(c-c^{-1}).$ Hence, by Corollary \ref{Qkl-cor}, for any $A_1,A_2\in M(4,\C),$
\begin{dmath}
Q_{k,l}(\imath_{c_1}(A_1),\imath_{c_2}(A_2))=k\cdot Q_{k-1,l}(\imath_{c_1}(A_1),\imath_{c_2}(A_2))\cdot i(c_1-c_1^{-1})+l\cdot Q_{k,l-1}(\imath_{c_2}(A_1),\imath_{c_2}(A_2))\cdot i(c_2-c_2^{-1}),
\end{dmath}
where $\imath_{c_i}$ are embeddings of $M(4,\C)$ into $M(2(k+l),\C)$ on the left and into $M(2(k+l)-2,\C)$ on the right.

In particular, 
$$Q_k(\imath_{c}(A))=i(c-c^{-1})\cdot kQ_{k-1}(\imath_{c}(A)),$$ 
implying part (1).

Similarly, by the Corollary \ref{Qkl-cor},
\begin{dmath}
Q_{k,1}(\imath_{c_1}(A_1),\imath_{c_2}(A_2))=
i(c_1-c_1^{-1}) k Q_{k-1,1}(\imath_{c_1}(A_1),\imath_{c_2}(A_2))+
i(c_2-c_2^{-1}) Q_k(\imath_{c_1}(A_1))
\end{dmath}
which by part (1) equals to 
$$i(c_1-c_1^{-1})k Q_{k-1,1}(\imath_{c_1}(A_1),\imath_{c_2}(A_2))+
\frac{1}{2}i^{k-1}(c_2-c_2^{-1})(c_1-c_1^{-1})^{k-2}k! Q_2(A_1).$$
Now we obtain part (2) of the statement by induction on $n$. 
\epr

\noindent{\it Proof of Proposition \ref{alpha-epi}:}  By \cite[Thm. 2]{S-quo}, $\C[X_{SO(4,\C)}(F_2)]$ is generated by trace functions and by the functions 
$$Q(\gamma,\gamma),\ \text{for}\ \gamma\in \{\gamma_1,\gamma_2,\gamma_1\gamma_2,\gamma_1\gamma_2^{-1}\}$$ 
and by 
\beq \lb{e_3Qs}
Q(\gamma_1,\gamma_2),\quad Q(\gamma_1,\gamma_2\gamma_1^{-1}),\quad 
Q(\gamma_2,\gamma_1\gamma_2^{-1}),
\eeq 
where $\gamma_1,\gamma_2$ generate $F_2$, as before. 
(We used here the fact that $Q$ is symmetric in its arguments.) 

Since 
$$\alpha_{c_1,c_2*}(\tau_\gamma)=\tau_\gamma+(c+c^{-1})(n-2),$$ 
for $c=c_1^{w_1(\gamma)}c_2^{w_2(\gamma)}$, 
the function $\alpha_{c_1,c_2*}$ maps the $SO(2n,\C)$-trace algebra of $F_2$ onto the $SO(4,\C)$-trace algebra. Therefore, it is enough to show that all of the above $Q$ functions belong to the image of $\alpha_{c_1,c_2*}$ as well.

Note that 
$$\alpha_{c_1,c_2*}Q_n(\gamma)[\rho]=Q_n((\imath_{c_1,c_2}\rho)(\gamma)).$$ 
By (\ref{e_alphac1c2}) and by Lemma \ref{Qimath}(1), that equals to 
$$\frac{1}{2}[i(c-c^{-1})]^{n-2}n! Q_2(\rho(\gamma))$$ 
for $c=c_1^{w_1(\gamma)}c_2^{w_2(\gamma)}.$
Therefore, up to a constant, $\alpha_{c_1,c_2*}Q_n(\gamma)$ equals $Q(\gamma,\gamma)$.
By the assumptions about $c_1$ and $c_2$, that constant is non-zero for 
$\gamma\in \{\gamma_1,\gamma_2,\gamma_1\gamma_2,\gamma_1\gamma_2^{-1}\}$.
Therefore, $Q(\gamma,\gamma)$ belongs to the image of $\alpha_{c_1,c_2*}$ for $\gamma$ as above.

Finally, it remains to be shown that the three functions of (\ref{e_3Qs}) belong to the image of $\alpha_{c_1,c_2*}$ as well. 

By (\ref{e_alphac1c2}) and by Lemma \ref{Qimath}(2), 
$$\alpha_{c_1,c_2*}Q_{n-1,1}(\gamma,\gamma')[\rho]=Q_{n-1,1}((\imath_{c_1,c_2}\rho)(\gamma),(\imath_{c_1,c_2}\rho)(\gamma'))$$ 
equals to 
$$\frac{1}{2}[i(c-c^{-1})]^{n-2}(n-1)! Q(\rho(\gamma),\rho(\gamma'))+d\cdot  Q(\rho(\gamma),\rho(\gamma)),$$
for $c=c_1^{w_1(\gamma)}c_2^{w_2(\gamma)}$. (Note that this ``$c$'' is the ``$c_1$'' in Lemma \ref{Qimath}(2) here.)
Consequently, 
$$\alpha_{c_1,c_2*}Q_{n-1,1}(\gamma,\gamma')=\frac{1}{2}[i(c-c^{-1})]^{n-2}n! Q(\gamma,\gamma')+d\cdot  Q(\gamma,\gamma).$$
By taking $$(\gamma,\gamma')\in \{(\gamma_1,\gamma_2),(\gamma_1,\gamma_2\gamma_1^{-1}), (\gamma_2,\gamma_1\gamma_2^{-1})\}$$
we see that linear combinations of the desired elements of (\ref{e_3Qs}) (with non-zero coefficients) and of $Q(\gamma_1,\gamma_1), Q(\gamma_2,\gamma_2)$ belong to the image of $\alpha_{c_1,c_2*}.$ Since we proved that $Q(\gamma_1,\gamma_1)$ and $Q(\gamma_2,\gamma_2)$ belong to the image of $\alpha_{c_1,c_2*}$ already, the proof is complete.
\qed

\bthm \lb{alpha-FT-onto}
For any group $\Gamma$ and any $c_1,c_2\in \C^*,$
$\alpha_{c_1,c_2*}$ maps ${\cal FT}_{SO(2n,\C)}(\Gamma)$ to ${\cal FT}_{SO(4,\C)}(\Gamma)$.
\ethm

For the sake of the proof of the above theorem it will be convenient to consider another matrix realization of $SO(2n,\C)$: Let $J_{2n}$ be a $2n\times 2n$ matrix composed of $n$ diagonal blocks of the form
$$J_2=\left(\begin{array}{cc} 0 & 1\\ 1& 0\\ \end{array}\right)$$
and let 
$$SO_J(2n,\C)=\{A: AJ_{2n}A^T=J_{2n},\ det(A)=1\}.$$  
Since $J_{2n}=K_{2n}\cdot K_{2n}^T,$ where $K_{2n}$ is composed of $n$ diagonal blocks of the form
$$K_2=\frac{1}{\sqrt{2}}\left(\begin{array}{cc} 1 & i\\ 1& -i\\ \end{array}\right),$$
the isomorphism 
\beq \lb{e_PhiSO_JSO}
\Phi(A)=K_{2n}^{-1}AK_{2n},\quad \Phi: SO_J(2n,\C)\to SO(2n,\C)
\eeq 
 identifies these two groups.\\

\noindent{\bf Proof of Theorem \ref{alpha-FT-onto}:}  It is enough to show that for every $\gamma\in \Gamma$, for every representation $\phi: SO(2n,\C)\to GL(N,\C),$ and for every $c_1,c_2\in \C^*$, $\alpha_{c_1,c_2*}(\tau_{\phi,\gamma})$ can be expressed by a polynomial in variables $\tau_{\psi_i,\gamma}$ for some $SO(4,\C)$-representations $\psi_i$. We will show (a seemingly stronger, but equivalent statement) that under the above assumptions,
\beq \lb{e_alphatau}
\alpha_{c_1,c_2*}(\tau_{\phi,\gamma})=\sum_{i=1}^k a_i \tau_{\psi_i,\gamma},
\eeq 
for some $SO(4,\C)$-representations $\psi_1,...,\psi_k$ and some $a_1,...,a_k\in \C.$

The map $\alpha_{c_1,c_2*}$ sends a regular function $f$ on $X_{SO(2n,\C)}(\Gamma)$  to a function on $X_{SO(4,\C)}(\Gamma),$ assigning to the equivalence class of $\rho:\Gamma\to SO(4,\C)$ the value $f(\alpha_{c_1,c_2}(\rho)).$ Therefore, (\ref{e_alphatau}) can be restated 
as
\beq \lb{e_alphatau2}
tr\phi(\alpha_{c_1,c_2}\rho)(\gamma)= \sum_{i=1}^k a_i tr(\psi_i(\rho(\gamma))).
\eeq

From now on we will identify $SO(2n,\C)$ with $SO_J(2n,\C)$ via (\ref{e_PhiSO_JSO}) and, hence, assume that $\rho:\Gamma\to SO_J(4,\C).$

Let $T_n$ be the group of diagonal matrices in $SO_J(2n,\C)$
such that the $(2i-1,2i-1)$ entry is the inverse of the $(2i,2i)$ entry
for $i=1,...,n.$ It is a maximal torus in $SO_J(2n,\C)$ (and its Cartan subgroup).
Since equation (\ref{e_alphatau2}) depends on the equivalence class of $\rho$ in $Hom(\Gamma,SO_J(4,\C))$ only and since this class contains a representative $\rho'$ such that $\rho'(\gamma)$ belongs to $T_2\subset SO_J(4,\C)$, we can assume for simplicity that $\rho(\gamma)$ belongs to $T_2$. 

Through the isomorphism (\ref{e_PhiSO_JSO}), 
the embedding $\imath_c: SO(4,\C)\to SO(2n,\C)$ corresponds to $\imath_c: SO_J(4,\C)\to SO_J(2n,\C)$ sending $A$ 
to a $2n\times 2n$ matrix composed of $A$ followed by $(n-2)$ diagonal blocks of the form 
$$\left(\begin{array}{cc} c & 0\\
0 & c^{-1}\\ \end{array}\right).$$
As in (\ref{e_alphac1c2}),
$$\alpha_{c_1,c_2}(\rho)(\gamma)=\imath_c\rho(\gamma),$$ 
where $c=c_1^{w_1(\gamma)}c_2^{w_2(\gamma)}$.
Hence, $\alpha_{c_1,c_2}(\rho)(\gamma)$ belongs to $T_n$.

Therefore, we can complete the proof of Theorem \ref{alpha-FT-onto} by taking $x=\rho(\gamma)\in T_2$ and by showing the following:

\blem \lb{l_alphatau3} 
For every representation $\phi$ of $SO(2n,\C)$ and every $c\in \C^*$, there are
representations $\psi_1,...,\psi_k$ of $SO(4,\C)$ and $a_1,...,a_k\in \C$ such that
\beq \lb{e_alphatau3}
tr\phi(\imath_c(x))= \sum_{i=1}^k a_i tr(\psi_i(x)),
\eeq
for all $x\in T_2.$
\elem 

\noindent{\it Proof:} Note that the left hand side of the above equation is a polynomial function on $T_2$ and the right hand side is a formal character of $SO(4,\C)$ with complex coefficients.

Note also that the weight lattice, $\Lambda_2,$ of $SO(4,\C)$ embeds into $\C[T_2]$ and that this embedding extends to $\C\Lambda_2 \to \C[T_2]$ which is an embedding as well, cf. \cite[III.8]{Bo}. The image of this map is an algebra generated by functions sending a diagonal matrix with diagonal entries $d_1,d_1^{-1},d_2,d_2^{-1}$ to $d_1^{\pm 2},$
$d_2^{\pm 2},$ and $d_1d_2$. Therefore, the image of $\C\Lambda_2$ is of index $2$ in 
$\C[T_2]$.

For complex reductive groups, the algebra of formal characters (over $\C$) coincides with the algebra $(\C\Lambda_2)^{W_2}$ of linear combinations of weights invariant under the Weyl group action, cf. \cite[\S 23.2]{FH}.

Therefore, the statement of Lemma \ref{l_alphatau3} follows from the following lemma.

\blem 
(1) $tr\phi(\imath_c(x))$ is $W_2$-invariant function of $x\in T_2.$\\
(2) $tr\phi(\imath_c(x))$ is a polynomial function on $T_2$ which belongs to $\C\Lambda_2.$
\elem 

\bpr (1) The Lie algebra $\frak t_n$ of $T_n$ has a basis $H_1,...,H_n,$ 
where $H_i$ is a $2n\times 2n$ matrix whose all entries are zero except for $(2i-1,2i-1)$-th entry $1$ and $(2i,2i)$-th entry $-1$. 
Since the exponential map $exp: \frak t_2\to T_2$ is onto, $x=exp(z),$ for some $z\in \frak t_2$. Then $\imath_c x=e^{z+dv}$ for
$v=H_3+...+H_n$ and for any $d$ such that $e^d=c.$

Suppose that $\phi$ has weights $\alpha_1,...,\alpha_N\in \Lambda_n$ with multiplicities $m_1,...,m_N.$ Then
$$tr\phi(\imath_c(x))=tr\phi(e^{z+dv})=\sum_{i=1}^N m_ie^{\alpha_i(z+dv)},$$
by \cite[(23.40)]{FH}.

Since the exponential map $exp: \frak t_2\to T_2$ is $W_2$-equivariant and onto, it is enough to show that the above expression is a $W_2$ invariant function of $z\in \frak t_2.$

Since $\alpha_i$ are weights of an $SO(2n,\C)$-representation and $m_i$ are their multiplicities,
$\sum_{i=1}^N m_ie^{\alpha_i(\cdot)}$ is a $W_n$-invariant function on $\frak t_n$. 
Hence, it is a linear combination of orbits of $W_n$,  i.e. a linear combination of expressions $\sum_{w\in W_n}  e^{w\cdot \alpha(\cdot)},$ for some $\alpha\in \frak t_n^*.$
Therefore, it is enough to prove that for every $\alpha\in \frak t_n^*$ and for every $d\in \C$,
$\sum_{w\in W_n}  e^{w\cdot \alpha(z+dv)}$ is $W_2$ invariant function of $z\in \frak t_2.$

Let $L_1,...,L_n\in \frak t_n^*$ be the dual basis to $H_1,...,H_n$.
Assume that $\alpha=\sum_{i=1}^n d_iL_i$, for some $d_1,...,d_n\in \C.$ 
Since $\alpha$ is a weight of an $SO(2n,\C)$-representation, it belongs to the weight lattice, implying that either $d_1,...,d_n\in \Z$ or $d_1,...,d_n\in \Z+\frac{1}{2}$.
Since elements of $W_n$ are signed permutations with even numbers of sign changes, elements of $W_n\alpha$ are of the form
$$\ve_1d_1L_{\sigma(1)}+...+\ve_nd_nL_{\sigma(n)},$$ 
for some $\sigma\in S_n$ and 
$\ve_1,...,\ve_n\in \{\pm 1\}$ such that $\ve_ i=-1$ for an even number of indices $i$.
By relabeling the indices by $\sigma^{-1}$, this expression can be written as 
$$\ve_1'd_{\tau(1)}L_1+...+\ve_n' d_{\tau(n)}L_n,$$ 
where $\tau=\sigma^{-1}$ and $\ve_ i'=\ve_{\tau(i)}$ for simplicity.

Since $v=H_3+...+H_n$ and $L_i(H_j)=\delta_{ij},$
\beq \lb{e_orbW}
\sum_{w\in W_n}  e^{w\cdot \alpha(z+dv)}=
\sum  e^{\ve_1'd_{\tau(1)}L_1(z)+\ve_2'd_{\tau(2)}L_2(z)+
d(\ve_3'd_{\tau(3)}+...+\ve_n'd_{\tau(n)})},
\eeq
where the sum on the right is over all $\tau\in S_n$ and all $\ve_1',...,\ve_n'\in \{\pm 1\}$ such that $\ve_ i'=-1$ for an even number of indices $i$.

Denote by $S_{k,l,r}(d)$ the sum
$$\sum  e^{d(\epsilon_1d_{\nu(1)}+...+\epsilon_{n-2}d_{\nu(n-2)})},$$
over all bijections $\nu: \{1,...,n-2\}\to \{1,...,n\}-\{k,l\},$ and all $\epsilon_1,...,\epsilon_{n-2}\in \{\pm 1\},$
such that $\#\{i: \ve_ i=-1\}=r$ mod $2$. Then by grouping the terms of (\ref{e_orbW}) by the values of $\tau(1)$ and of $\tau(2)$, the sum (\ref{e_orbW}) can be written as
\beq \lb{e_orbW2}
\sum_{k\ne l} \left(\sum_{\epsilon=\pm 1} e^{\epsilon d_{k}L_1(z)+\epsilon d_{l} L_2(z)} S_{k,l,0}(d)+
\sum_{\epsilon=\pm 1} e^{\epsilon d_{k}L_1(z)-\epsilon d_{l} L_2(z)} S_{k,l,1}(d)\right).
\eeq 

Since we can choose two generators of $W_2=\Z_2\times \Z_2$ such that the first one applied to $z$ interchanges the values of $L_1(z)$ and of $L_2(z)$ and the second negates $L_1(z)$ and $L_2(z)$, it is easy to see that  
\beq \lb{e_orbW3}\sum_{\epsilon=\pm 1} e^{\epsilon d_{k}L_1(z)+\epsilon d_{l} L_2(z)} \quad
\text{and}\quad  
\sum_{\epsilon=\pm 1} e^{\epsilon d_{k}L_1(z)-\epsilon d_{l} L_2(z)}
\eeq
are $W_2$-invariant functions of $z$.
Consequently, (\ref{e_orbW}) is a $W_2$-invariant function of $z$ and the first statement of lemma follows.

(2) By the above argument $tr\phi(\imath_c(x))$ can be written as a linear combination of expressions (\ref{e_orbW3}), for $d_k,d_l\in \Z$ or $d_k,d_l\in \Z+\frac{1}{2}.$ Since $L_1,L_2, \frac{1}{2}(L_1+L_2)$ are weights of $\frak t_2$, the above expressions belong to $\C[\Lambda_2]\subset \C[T_2].$
\epr

This completes the proof of Lemma \ref{l_alphatau3} and of Theorem \ref{alpha-FT-onto}.
\qed

 \section{Proof of Theorem \ref{undisting}}

Consider the symmetric square, $Sym^2(\C^5)$, which is a $15$-dimensional space composed of vectors 
$$v\odot w=v\otimes w+w\otimes v\in \C^5\otimes \C^5.$$ 
The standard inner product on $\C^5$ yields an inner product on $\C^5\otimes \C^5$,
$$(v_1\otimes v_2,w_1\otimes w_2)=(v_1,w_1)(v_2,w_2)$$ 
which restricts to an inner product on $Sym(\C^5).$ If $e_1,...,e_5$ is the standard orthonormal basis of $\C^5$ then $e_i\odot e_j,$ for $i\leq j,$ is a basis of $Sym^2(\C^5)$
and an easy computation shows that 
\beq \lb{e_scalar}
(e_i\odot e_j,e_k\odot e_l)=\begin{cases} 4\ \text{for}\ i=k=l=j\\
2\ \text{for}\ i=k\ne j=l\\
0\ \text{otherwise.}\\
\end{cases}
\eeq 
Note that the natural action of $SO(5,\C)$ on $Sym^2(\C^5)$ preserves the above inner product. 

Let $e^1,...,e^5$ the dual basis of $(\C^5)^*.$
Then $tr: \C^5\otimes (\C^5)^*\to \C$  which can be represented by the tensor
$$\sum_{i=1}^5 e^i\otimes e_i\in (\C^5\otimes (\C^5)^*)^*=(\C^5)^*\otimes \C^5$$
is $SO(5,\C)$ invariant. Since the map from $\C^5$ to $(\C^5)^*$ sending $e_i$ to $e^i$ is a $SO(5,\C)$-module isomorphism, we can identify this tensor with 
$$z=\sum_{i=1}^5 e_{i}\otimes e_i.$$
Therefore, as an $SO(5,\C)$-module, $Sym^2(\C^5)$ splits into $\C z$ and $z^{\perp}$. By \cite[Ex 24.32]{FH}, the representation $SO(5,\C)\to SO(z^{\perp})=SO(14,\C)$
is irreducible. We will denote that representation by $\alpha.$ 

For the proof of Theorem \ref{undisting} it is enough to assume that $\Gamma=\Z_p*\Z_q.$

By Corollary \ref{phipsi-irred} below there is a representation $\psi: \Gamma\to SO(5,\C)$ such that $\rho=\alpha\psi:\Gamma\to SO(14,\C)$ is irreducible. 
For $n=7$, let $\rho=\alpha\psi$. If $n\geq 9$ then let $\rho=\alpha\psi\oplus \eta$, where $\eta: \Z_p*\Z_q\to SO(2n-14)$ is an irreducible representation, whose existence is implied by Corollary \ref{eta-irred}.

Let $\rho'$ be obtained from $\rho$ by the involution $\sigma$ of \cite{S-char}, that is by a conjugation by a matrix $M\in O(14,\C)$ of determinant $-1$.
Clearly, $\rho$ and $\rho'$ are undistinguishable by trace functions.  
Furthermore, we claim that for $\rho$ as above $Q_n(\gamma)(\rho)$ (that is $Q_n(\rho(\gamma))$) vanishes for all $\gamma\in \Gamma.$ 
To see that, note that $\psi(\gamma)$, being an element of $SO(5,\C),$ has eigenvalues $1,x^{\pm 1},y^{\pm 1},$ for some $x,y\in \C^*.$ The eigenvalues of the action of $\psi(\gamma)$ on $Sym(\C^5)$ are products of the above ones and, in particular, the eigenvalue $1$ appears with multiplicity at least three -- coming from 
$1\cdot 1, x\cdot x^{-1}$, and from $y\cdot y^{-1}.$ Consequently, the eigenspace of $1$ for the action of $\psi(\gamma)$ on $z^\perp$ has dimension at least $2.$ Now, by Lemma \ref{prop_Q} and the fact that different eigenspaces are orthogonal to each other and that
$Q_2\left(\begin{array}{cc} 1 & 0\\ 0& 1\\ \end{array}\right)=0,$ we have $Q_n(\gamma)=0$ indeed for all $\gamma\in \Gamma.$

Now, it remains to be shown that $\rho$ and $\rho'$ are not equivalent. Since they are both completely reducible, it is enough to show that they are not conjugate in $SO(2n,\C).$ (By \cite[Thm. 30]{S-char}, completely reducible representations have closed conjugacy orbits in the representation variety and, therefore, are non-equivalent in the character variety if they are not conjugate.) Suppose $\rho'=M\rho M^{-1}$ equals $A\rho A^{-1},$ for some $A\in SO(2n,\C)$. Then  $MA^{-1}$ commutes with the image of $\rho,$ which by Shur's Lemma spans the entire matrix algebra $M(14,\C)$ if $n=7$ or, spans $M(14,\C)\oplus M(2n-14,\C)$, for $n\geq 9.$ Hence, for $n=7$, $MA^{-1}$ is a scalar matrix. The only scalar matrices in $O(2n,\C)$ are $\pm I.$ That however implies that $det(MA^{-1})=1$ and, hence, 
$$det(A)=det(M)= -1,$$
contradicting the assumption that $\rho$ and $\rho'$ are conjugate in $SO(2n,\C).$ 

For $n\geq 9,$ $MA^{-1}$ is composed of a $14\times 14$ scalar matrix and of another $(2n-14)\times (2n-14)$ scalar matrix. The only such matrices in $O(2n,\C)$ are $(\pm I_{14})\oplus (\pm I_{2n-14})$. Their determinants are $1$, implying that $det(A)= -1$ again.\vspace*{.1in}
\qed

The above proof relies on two statements about irreducible representations of $\Z_p*\Z_q$. We will start with the proof of the easier one, asserting the existence of irreducible representations of $\Z_p*\Z_q$ into $SO(2m,\C)$ for any $m>1.$ The method of proof of this statement will be extended to show the second statement (Lemma \ref{ZpZq} and Corollary \ref{phipsi-irred}) afterwards.

Let $B_1$ be a block matrix composed of $D_{\xi_p},$ $D_{\xi_p^2},$ ..., $D_{\xi_p^m},$
where $\xi_p$ is a primitive $p$-th root of $1$,
and, similarly, let $B_2$ be a block matrix composed of $D_{\xi_q},$ $D_{\xi_q^2},$ ...,$D_{\xi_q^m},$
where $\xi_q$ is a primitive $q$-th root of $1$.
Since $B_1,$ $B_2$ are of order $p$ and $q$ respectively, for every 
$A\in SO(2m,\C)$
there is a representation 
$$\eta_A: \Z_p*\Z_q\to SO(2m,\C)$$ 
sending a generator of $\Z_p$ to $B_1$ and a generator of $\Z_q$ to $AB_2A^{-1}.$
Assume that $p,q> 2m.$ Since $B_1,$ $B_2$ commute and have $2m$ distinct eigenvalues (each), they share $2m$ distinct eigenvectors $v_1,...,v_{2m}.$ 
Since $B_1,B_2$ are orthogonal transformations, $v_1,...,v_{2m}$ form an orthogonal basis of $\C^{2m}.$

Given $v,w\subset \C^{2m}$, let $Z_{v,w}$  denote the set of matrices $A\in SO(2m,\C)$ such that $v\subset A w^{\perp}.$

\bpro \lb{ZpZq-irred} 
(1) For every $m>2$ and $p,q> 2m,$ if $\eta_A: \Z_p*\Z_q\to SO(2m,\C)$ is reducible then $A\in \bigcup_{k,l=1}^{2m} Z_{v_k,v_l}.$\\
(2) For every non-zero $v,w\in \C^{2m}$, $Z_{v,w}$ is Zariski closed proper subspace of $SO(2m,\C).$
\epro

\bcor \lb{eta-irred}
For every $m>2$ and $p,q> 2m,$ $\eta_A: \Z_p*\Z_q\to SO(2m,\C)$ is irreducible 
for non-empty Zariski open set of matrices $A$ in $SO(2m,\C).$ 
\ecor 

\noindent {\it Proof of Proposition \ref{ZpZq-irred}:}\\
(1) Suppose that $\eta_A$ is reducible. Then $B_1$ and $AB_2A^{-1}$ preserve a certain proper non-empty subspace $V\subset \C^{2m}.$  Since $V$ is preserved by $B_1$, it must have a basis given by a subset of $\{v_1,...,v_{2m}\}.$ Suppose then that $v_k\in V.$ Since $AB_2A^{-1}$ preserves $V$, the transformation
$B_2$ preserves $A^{-1}V$ implying that $A^{-1}V$ has a basis which is also a subset of $\{v_1,...,v_{2m}\}$ (the eigenvectors of $B_2$). Since this basis is a proper subset of $\{v_1,...,v_{2m}\}$, $v_l\not\in A^{-1}V$ for some $l.$ Since all $v_i,$ for $i\ne l$, are orthogonal to $v_l,$ we have $A^{-1}v_k \perp v_l$ and, hence, 
$A\in Z_{v_k,v_l}.$\\
(2) It is clear that $Z_{v,w}$ is Zariski closed. Since $SO(2m,\C)$ acts transitively on $\C^{2m},$ the set $Z_{v,w}$ is a proper subset of $SO(2m,\C)$.\vspace*{.1in}
\qed

Recall that $z^\perp$ is a $14$-dimensional subspace of $Sym^2(\C^5)$ and that
$SO(5,\C)$ acts on it through an irreducible representation $\alpha.$
The reminder of this section is devoted to proving that there is a representation $\psi: \Gamma\to SO(5,\C)$ such that $\rho=\alpha\psi:\Gamma\to SO(z^\perp)$ is irreducible.  (That is the last component of the proof of Theorem \ref{undisting} which needs to be established).

Let $B_c\in SO(5,\C)$ for $c\in \C^*$ be a block matrix composed of $D_{c},$ $D_{c^4},$ and of the $1\times 1$ identity matrix along the diagonal. (Note that the meaning of $B$ here is different, but analogous to that of the previous proof.)

Let $\psi_A: \Z_p*\Z_q\to SO(5,\C)$ be a representation sending a generator of $\Z_p$ to $B_{\xi_p}$ and a generator of $\Z_q$ to $AB_{\xi_q}A^{-1},$ for some $A\in SO(5,\C),$ where $\xi_p$ and $\xi_q$ are primitive $p$-th and $q$-th roots of $1$ respectively, as before. (These matrices have order $p$ and $q$, respectively, and, therefore, $\psi_A$ is well defined.)

Since $B_{\xi_p}$ has eigenvalues $\xi_p^{\pm 1},$ $\xi_p^{\pm 4},$ and $1$, the eigenvalues of $B_{\xi_p}$ acting on $Sym^2(\C^5)$ are monomials of degree two in the above eigenvalues, that is  
\beq \lb{e_eiegnnot1}
\xi_p^{\pm 2}, \xi_p^{\pm 8}, \xi_p^{\pm 1 \pm 4}, \xi_p^{\pm 1}\cdot 1, \xi_p^{\pm 4}\cdot 1,
\eeq 
and 
$$\xi_p\cdot \xi_p^{-1}=\xi_p^4\cdot \xi_p^{-4}=1\cdot 1=1.$$
Assume now that $p,q>16.$
Since $\xi_p$ is a primitive $p$-th root of unity, all of these eigenvalues appear once, except for $1$ which appears with multiplicity $3$.
Let $v_1,...,v_{12}$ be eigenvectors for the eigenvalues (\ref{e_eiegnnot1}).
Denote the eigenspace of $1$ inside $z^\perp$ by $F.$  Then $\dim F=2$ and $v_1,...,v_{12},z$ and $F$ are all orthogonal to each other.

Note that $v_1,...,v_{12}$ are eigenvectors of the action of $B_{\xi_q}$ on $Sym^2(\C^5)$ as well and that $F$ is also the eigenspace of $1$ for the $B_{\xi_q}$-action on $z^\perp.$

Denote by ${\cal R_{p,q}}\subset SO(5,\C)$ the set of matrices $A$ such that the representation $\alpha\psi_A: \Z_p*\Z_q\to SO(z^\perp)$ is reducible. 

Let $Z_{v,w}$, for $v,w\subset z^\perp,$ denote the set of matrices $A\in SO(5,\C)$ such that $\alpha(A)v\subset w^{\perp}.$
Let  $Y$ denote the set of matrices $A\in SO(5,\C)$ such that the dimension of the subspace $F+\alpha(A)F$ of $z^\perp$ is at most $3$.

\blem \lb{ZpZq} 
Under the above assumptions\\
(1) For every nonzero $v,w\in \C^5$, $Z_{v,w}$ is a Zariski closed proper subset of $SO(5,\C).$\\
(2) $Y$ is a Zariski closed proper subset of $SO(5,\C).$\\
(3) 
$${\cal R_{p,q}}\subset \bigcup_{1\leq i,j\leq 12} Z_{v_i,v_j}\cup Y.$$
\elem

\bcor \lb{phipsi-irred}
${\cal R_{p,q}}$ is a subset of a proper closed subset of $SO(5,\C)$ and, hence, $\alpha\psi_A$ is irreducible for a non-empty Zariski-open set of matrices $A\in SO(5,\C).$
\ecor 

\noindent {\it Proof of Lemma \ref{ZpZq}:}
(1) Since $\alpha$ is irreducible, for every nonzero $v\in z^\perp$ the set $\alpha(SO(5,\C))v$ spans $ z^\perp$. Therefore, $Z_{v,w}$ is a proper subset of $SO(5,\C).$ Clearly, it is Zariski closed.

(2) Since $\begin{pmatrix}1\\1\end{pmatrix}$ and $\be{pmatrix}  1 \\ -1 \en{pmatrix}$ are eigenvectors of $D_c$ with eigenvalues $c$ and $c^{-1}$ respectively,
$$f_1=\begin{pmatrix}1\\1\\ 0\\ 0\\ 0\end{pmatrix}\odot \be{pmatrix}  1 \\ -1\\0 \\0\\ 0 \en{pmatrix}\quad \text{and}\quad f_2= \begin{pmatrix}0\\ 0\\ 1\\1\\ 0\end{pmatrix}\odot \be{pmatrix}  0\\0\\ 1 \\ -1\\0\en{pmatrix}$$ 
are eigenvectors with eigenvalue $1$ of $\alpha(D_{\xi_p})$ and of $\alpha(D_{\xi_q}).$
By (\ref{e_scalar}), these vectors are orthogonal to each other and one can also easily check that they are orthogonal to $z$ as well. Therefore, $f_1,f_2$ form an orthogonal basis of $F.$
The condition
$$\dim F+\alpha(A)F\leq 3$$
is equivalent to the rank of the $4\times 14$ matrix built of coordinates of $f_1,$ $f_2,$ $\alpha(A)(f_1),$ $\alpha(A)(f_2)$ (with respect to some basis of $z^\perp$)
being at most $3$. Since that corresponds to vanishing of one of a finite number of $4\times 4$ minors of that matrix, $Y$ is algebraically closed. Therefore, it is remains to be shown that $A\not\in Y$ for some $A\in SO(5,\C)$. 
We claim that the cyclic permutation matrix 
$$A=
\begin{pmatrix}
0 & 1 & 0 & 0 & 0\\
0 & 0 & 1 & 0 & 0\\
0 & 0 & 0 & 1 & 0\\
0 & 0 & 0 & 0 & 1\\
1 & 0 & 0 & 0 & 0\\
\end{pmatrix}$$
will do the job.

Note that the action of $A$ on $Sym^2(\C^5)$ sends $f_1$ and $f_2$ to
$$g_1=\begin{pmatrix}0\\ 1\\1\\ 0\\ 0\end{pmatrix}\odot \be{pmatrix}  0\\ 1 \\ -1\\0 \\0 \en{pmatrix}\quad \text{and}\quad
g_2= \begin{pmatrix}0\\ 0\\ 0\\ 1\\1\end{pmatrix}\odot \be{pmatrix}  0\\0\\0\\ 1 \\ -1\en{pmatrix},$$
respectively. 
Since $f_1,f_2,g_1,g_2$ are linearly independent and orthogonal to $z$,
$$\dim F+\alpha(A)F=4.$$ 

(3) Let $A\in {\cal R}_{p,q}.$ Then $\alpha\psi_A$ is reducible and, hence,
$\alpha(B_{\xi_p})$ and $\alpha(AB_{\xi_q}A^{-1})$ preserve some non-zero, proper subspace $V$ of $z^\perp$ and its orthogonal complement $V^{\perp}$ in $z^\perp.$ Since $V$ is preserved by $\alpha(B_{\xi_p})$, it is a sum of a subspace of $F$ and of a space with a basis being a subset of $\{v_1,...,v_{12}\}.$

By replacing $V$ by $V^\perp$ if necessary, we can assume without loss of generality 
that $\dim V\geq \dim V^\perp,$ i.e. $\dim V\geq 7.$
Therefore, $V$ contains $v_i$ for some $i=1,...,12.$ If $\alpha(A)^{-1}V$ does not contain some $v_j$ then $A^{-1}v_i\perp v_j$ and $A\in Z_{v_i,v_j},$ by an argument analogous to that in Proposition \ref{ZpZq-irred}(1). Therefore, assume now that $\alpha(A)V$ contains all $v_i$, $i=1,...,12.$ 

Furthermore, if $V$ does not contain some $v_i$ then $v_i\perp V$ and, consequently,
$A$ sends it to a vector orthogonal to $AV$. Hence, by the above argument,
$A\in Z_{v_i,v_j},$ for every $j.$

Therefore, we can assume that both $V$ and $\alpha(A)V$ contain $v_1,...,v_{12}.$ Then $V^\perp$ is either $F$ or a $1$-dimensional subspace of $F$ and $\alpha(A)$ maps it to some subspace of $F.$
That implies that 
$$\dim F+\alpha(A)F\leq 3.$$
(Note that $F+\alpha(A)F$ is indeed three-dimensional if $\dim V^\perp=1$ and $\alpha(A)$ maps a complementary space to it in $F$ to outside of $F$.)
\qed

 \section{Proof of Proposition \ref{normalization}}

\bpr (1) By \cite[Prop 18]{S-quo}, $\cal FT_{SO(2n,\C)}(\Gamma)$ is generated  by $Q(\gamma,...,\gamma)$ for $\gamma\in \Gamma$ as an algebra over $\cal T_{SO(2n,\C)}(\Gamma)$. Therefore, by the assumption of the proposition, 
\beq \lb{e_QnotinT}
Q(\gamma,...,\gamma)\not\in \cal T_{SO(2n,\C)}(\Gamma)
\eeq 
for some $\gamma.$ By \cite{S-quo}, 
\beq \lb{e_QQT}
Q(\gamma_1,...,\gamma_n)Q(\gamma,...,\gamma)\in \cal T_{SO(2n,\C)}(\Gamma)
\eeq 
for any $\gamma_1,...,\gamma_n\in \Gamma$ implying that there is an embedding
$$\alpha: \C[X_{SO(2n,\C)}(\Gamma)]\to \frac{1}{Q(\gamma,...,\gamma)}\cal T_{SO(2n,\C)}(\Gamma),$$
defined on generators of $\C[X_{SO(2n,\C)}(\Gamma)]$ by
$$\alpha(\tau_\gamma)=\tau_\gamma,\quad \alpha(Q(\gamma_1,...,\gamma_n))=\frac{Q(\gamma_1,...,\gamma_n)Q(\gamma,...,\gamma)}{Q(\gamma,...,\gamma)}$$
which extends to an isomorphism between the field of fractions of $\C[X_{SO(2n,\C)}(\Gamma)]$ and the field of fractions of $\cal T_{SO(2n,\C)}(\Gamma).$ This proves the birationality of $\psi.$ 
The finiteness of $\psi$ follows from the fact that the elements $Q(\gamma_1,...,\gamma_n)$, for $\gamma_1,...,\gamma_n\in \Gamma$, which generate $\C[X_{SO(2n,\C)}(\Gamma)]$ as an $\cal FT_{SO(2n,\C)}(\Gamma)$-algebra, satisfy 
$$Q(\gamma_1,...,\gamma_n)^2\in \cal T_{SO(2n,\C)}(\Gamma),$$ 
by Corollary 14 and Corollary 17 of \cite{S-quo}, and, therefore,  are integral over $\cal T_{SO(2n,\C)}(\Gamma)$. 

Every finite, birational map is generically $1$-$1$ and it is a normalization map if its domain is normal.

(2) Since the involution $\sigma$ of \cite{S-quo} fixes $\cal T_{SO(2n,\C)}(\Gamma)$ and negates every $Q(\gamma,...,\gamma),$ condition (\ref{e_QnotinT}) is equivalent to
$$Q(\gamma,...,\gamma)\ne 0.$$
Since 
$$Q_n(\gamma)=(2i)^n n! (\tau_{D^+_n,\gamma}-\tau_{D^-_n,\gamma})$$
by \cite[Prop. 10]{S-gen}, and,
$$\tau_{D^+_n,\gamma}+\tau_{D^-_n,\gamma}\in \cal T_{SO(2n,\C)}(\Gamma),$$
we can substitute $Q_n(\gamma)$ for $\tau_{D^\pm_n,\gamma}$ in the statement of the proposition without loss of generality.
Since the elements of $\C[X_{SO(2n,\C)}(\Gamma)]$ distinguish non-equivalent representations and $\C[X_{SO(2n,\C)}(\Gamma)]$ is generated by functions $Q(\gamma_1,...,\gamma_n)$ as an $\cal T_{SO(2n,\C)}(\Gamma)$-algebra and since 
$$Q(\gamma_1,...,\gamma_n)([\rho])=\frac{f([\rho])}{Q(\gamma,...,\gamma)([\rho])},$$
for some $f\in \cal T_{SO(2n,\C)}(\Gamma)$, as long as $Q(\gamma,...,\gamma)([\rho])\ne 0,$
the trace functions together with $Q(\gamma,...,\gamma)$ distinguish $\rho$ from all non-equivalent representations as long as $Q(\gamma,...,\gamma)([\rho])$ does not vanish. Since
$Q(\gamma,...,\gamma)\ne 0$, that is the case for a generic $\rho$.
\epr

%

\end{document}